\documentclass[a4paper,12pt]{amsart}

\usepackage{xypic}

\newtheorem{theo}{Theorem}[section]
\newtheorem{ex}[theo]{Example}
\newtheorem{exs}[theo]{Examples}
\newtheorem{prop}[theo]{Proposition}
\newtheorem{lem}[theo]{Lemma}

\newtheorem{facts}[theo]{Facts}

\newtheorem{cor}[theo]{Corollary}

\newtheorem{rema}[theo]{Remark}

\def \Romannumeral #1 {\expandafter\uppercase\expandafter {\romannumeral #1} }

\def \Ga {{\Gamma}}

\def \gal {{\rm{Gal}}}

\def \calo {{\mathcal O}}
\def \spec {{\rm{Spec\,}}}

\def\lra{\longrightarrow}

\def \Z {{\bf Z}}

\def \C {{\bf C}}

\def \cok {{\rm{coker\,}}}

\def \G {{\bf G}_m}
\def \Ga {{\bf G}_a}
\def \A {{\bf A}}

\def\smallsquare{\vbox{\hrule\hbox{\vrule height 1 ex\kern 1 ex\vrule}\hrule}}
\def\enddem{\hfill \smallsquare\vskip 3mm}

\usepackage{amscd}
\usepackage{amssymb}
\usepackage{amsmath}

\everymath{\displaystyle}

\DeclareFontFamily{U}{wncy}{}
\DeclareFontShape{U}{wncy}{m}{n}{%
   <5>wncyr5%
   <6>wncyr6%
   <7>wncyr7%
   <8>wncyr8%
   <9>wncyr9%
   <10>wncyr10%
   <11>wncyr10%
   <12>wncyr6%
   <14>wncyr7%
   <17>wncyr8%
   <20>wncyr10%
   <25>wncyr10}{}
\DeclareMathAlphabet{\cyrille}{U}{wncy}{m}{n}

\xyoption{curve}

\def \sc {^{\textup{sc}}}
\def \ss {^{\textup{ss}}}
\def \ssu {^{\textup{ssu}}}

\def \u {^{\textup{u}}}
\def \red {^{\textup{red}}}
\def \ad {^{\textup{ad}}}
\def \lin {^{\textup{lin}}}
\def \sab {^{\textup{sab}}}
\def \ab {^{\textup{ab}}}
\def \SA {{\textup{SA}}}

\title{\'Etale homotopy groups of algebraic groups and homogeneous spaces}
\author{Cyril Demarche and Tam\'as Szamuely}

\address{Sorbonne Universit\'e and Universit\'e Paris Cit\'e, CNRS, IMJ-PRG, F-75005 Paris, France.}
\email{cyril.demarche@imj-prg.fr}

\address{Dipartimento di Matematica, Universit\`a di Pisa, Largo Bruno Pontecorvo 5, 56127 Pisa, Italy}
\email{tamas.szamuely@unipi.it}

\date{\today}
\begin{document}
\maketitle
\markright{\'ETALE HOMOTOPY GROUPS OF ALGEBRAIC GROUPS}

\begin{abstract}
We show the vanishing of the second homotopy group of the \'etale homotopy type of a smooth connected algebraic group over a separably closed field, completed away from the characteristic. This is an algebraic analogue of a classical theorem of Elie Cartan. Based on this result, we establish an explicit formula for the similarly completed second homotopy group of a homogeneous space.
\end{abstract}

\section{Introduction}

Computing the (unstable) homotopy groups of compact Lie groups is a much-studied classical problem. A lot is known on the subject but only low-degree results have some uniformity. The behaviour of higher degree groups is much more irregular, and they are computed case by case following the classification of simple Lie groups (see e.g. \cite{mimura}, \S 3.2  for a survey).

Perhaps the most famous uniform theorem is a classical result of \'Elie Cartan \cite{cartan}: all compact Lie groups have trivial second homotopy. Our first main result in this note is the following analogue in algebraic geometry.

\begin{theo}\label{pi2p} Let $G$ be a connected smooth algebraic group over a separably closed field $k$ of characteristic $p\geq 0$. Then $\pi_2(G^{\wedge (p')},1)=0$.
\end{theo}

Here and in what follows for a $k$-scheme $X$ the notation $ X^{\wedge (p')}$ stands for the completion of the \'etale homotopy type of $X$ with respect to the class of finite groups of order prime-to-$p$ as defined by Artin and Mazur \cite{artinmazur}. We then take the second homotopy group pointed at the unit element $1\in G(k)$. Note that since our $G$ is smooth, its \'etale homotopy type is profinite (see Fact \ref{ethomfacts} (1) below), so the completion operation does not change anything in characteristic 0, whereas it is of crucial importance in positive characteristic, as we shall see.

In characteristic 0 the theorem will be deduced from that of Cartan using comparison theorems between classical and \'etale homotopy. The proof of the positive characteristic case is more involved, however, and is done by specialization and fibration techniques.

One of the difficulties is caused by the fact that the fibration exact sequence in \'etale homotopy theory is known to exist only under rather restrictive assumptions (see the next section). The following theorem provides an important case where it applies:

\begin{theo}\label{homex} Let $k$ be an algebraically closed field, and $G$ a smooth connected algebraic group over $k$. Let $H\subset G$ be a closed connected subgroup, and denote by $X$ the quotient $G/H$. There is a long exact sequence
 $$\cdots \to  \pi_{i+1}( X^{\wedge (p')}, \bar 1) \to \pi_i( H^{\wedge (p')},1) \to \pi_i( G^{\wedge (p')},1) \to \pi_i( X^{\wedge (p')},\bar 1) \to \cdots $$
 of \'etale homotopy groups, where $\bar 1$ is the image of ${1\in G(k)}$ in $X(k)$.
\end{theo}

\begin{rema}\label{Hred}\rm
The theorem holds for a general closed connected subgroup scheme $H\subset G$ but we may assume $H$ is equipped with its reduced structure. Indeed, in the general case we may consider the reduced subgroup scheme $\widetilde H\subset H$ and the quotient map $G/\widetilde H\to G/H$. As this map is finite and purely inseparable, it induces an equivalence of \'etale sites, hence an isomorphism of (completed) \'etale homotopy types.
\end{rema}

Note that according to Theorem \ref{pi2p}, and assuming $H$ is smooth, the above exact sequence breaks up in two segments. One ends by
\begin{equation}\label{pi3up}\cdots \to  \pi_{4}( X^{\wedge (p')}, \bar 1) \to \pi_3( H^{\wedge (p')},1)\to \pi_3( G^{\wedge (p')},1)\to \pi_3( X^{\wedge (p')},\bar 1) \to 0, \end{equation}
the other is
\begin{equation}\label{pi2down}0 \to \pi_2( X^{\wedge (p')},\bar 1) \to \pi_1(H,1)^{(p')}\to \pi_1(G,1)^{(p')} \to \pi_1(X,\bar 1)^{(p')} \to 0. \end{equation}

Here for the last three terms we have used the fact (\cite{artinmazur}, Corollary 3.7) that for a smooth $k$-scheme $X$ and $i=1$ we have an isomorphism $\pi_1( X^{\wedge (p')})\cong  \pi_1(X)^{(p')}$ where the latter group is the maximal prime-to-$p$-quotient of the \'etale fundamental group.

Using sequence (\ref{pi2down}) we can compute $\pi_2( X^{\wedge (p')},\bar 1)$ more precisely. This is done by breaking up $G$ and $H$ in pieces. By Chevalley's theorem (see e.g. \cite{brionbook}, Chapter 2) the group $G$ has a maximal closed connected linear subgroup $G^{\rm lin}$; denote by $G^{\rm u}$ its unipotent radical over $k$. The derived subgroup $G^{\rm ss}$ of the pseudo-reductive quotient $G^{\rm lin}/G^{\rm u}$ has a simply connected cover $G^{\rm sc}$; denote by $T_{G^{\rm sc}*}$ the cocharacter group of a maximal torus $T_{G\sc}$ in $G\sc$.

On the other hand, once a maximal torus $T_G\subset G^{\rm lin}/G^{\rm u}$ is fixed containing the image of $T_{G\sc}$, it can be embedded in a unique maximal semi-abelian  variety ${\rm SA}_G$ contained in $G/G^{\rm u}$; it is an extension of the maximal abelian variety quotient $G^{\rm ab}$ of $G$ by $T_G$. Changing $T_G$ amounts to replacing ${\rm SA}_G$ by a conjugate subgroup. For a construction of  ${\rm SA}_G$, see (\cite{dem:brauer}, \S 4.1) or Remark \ref{gantrema} below.

Denote by $T_{(p')}({\rm SA}_G)$ the prime-to-$p$ Tate module of ${\rm SA}_G$. Recall that this profinite abelian group is the extension of the Tate module $T_{(p')}(G^{\rm ab})$ by the Tate module $T_{(p')}(T_G)$ and the latter group is just the cocharacter group ${T_G}_*$ of $T_G$ tensored by $\Z_{(p')}(1) $, the inverse limit of all prime-to-$p$ roots of unity in $k$. In particular, by our choice of $T_G$ we have a map $\tau_G:\,{T_{G\sc}}_* \otimes \Z_{(p')}(1) \to T_{(p')}(\SA_G)$. The following statement generalizes Proposition 3.10 of \cite{dem}.

\begin{prop}\label{proppi1G} Still assuming $k$ algebraically closed,
there exists a canonical short exact sequence of profinite abelian groups
$$
0 \to {T_{G\sc}}_* \otimes \Z_{(p')}(1)\stackrel{\tau_G} \lra T_{(p')}(\SA_G)\to \pi_1(G,1)^{(p')} \to 0.
$$
\end{prop}

Now apply the proposition for $G$ and $H$ (the latter assumed to be smooth thanks to Remark \ref{Hred}). Plugging the resulting expression in exact sequence  (\ref{pi2down}) gives an explicit description of $\pi_2( X^{\wedge (p')},\bar 1)$. Namely, introduce the complex of profinite abelian groups
 $$\mathcal{C}_{X,p'} := [{T_{H\sc}}_* \otimes \Z_{(p')}(1) \to T_{(p')}(\SA_H) \oplus ({T_{G\sc}}_* \otimes \Z_{(p')}(1) )\to T_{(p')}(\SA_G)]$$
placed in homological degrees 2, 1 and 0. The maps in the complex come from $\tau_G$ and $\tau_H$ defined above and from choosing maximal tori in $H\sc$ and $G\sc$ that are compatible via the map $H\sc\to G\sc$ induced by the universal property of the simply connected cover.

\begin{cor}\label{corpigp}  There is a canonical isomorphism of abelian profinite  groups
$$\pi_2( X^{\wedge (p')},1)\cong H_{1}(\mathcal{C}_{X,p'}).  $$

\end{cor}

\begin{rema}\rm Define $\mathcal{C}_{X\lin}$ as the three-term complex
$$\mathcal{C}_{X\lin} := [{T_{H\sc}}_* \to {T_H}_* \oplus {T_{G\sc}}_* \to {T_G}_*]$$
of free abelian groups of finite rank, with maps induced by suitable choices of maximal tori as explained above. When $G$ is linear, we plainly have an identification of complexes $$\mathcal{C}_{X\lin}\otimes \Z_{(p')}(1)\cong \mathcal{C}_{X,p'}.$$ However, if we only assume $H$ to be linear, we still have an isomorphism $H_{1}(\mathcal{C}_{X\lin})\otimes \Z_{(p')}(1) \stackrel\sim\to H_{1}(\mathcal{C}_{X,p'})$ due to the
injectivity of the natural map ${T_G}_*\otimes\Z_{(p')}(1)\to T_{(p')}(\SA_G)$.
Therefore Corollary \ref{corpigp} gives
\begin{equation}\label{clin}H_{1}(\mathcal{C}_{X\lin}) \otimes \Z_{(p')}(1) \xrightarrow{\sim} \pi_2(X^{\wedge (p')}).\end{equation}
The stabilizer $H$ is known to be linear when the action of $G$ is faithful, by a theorem of Matsumura \cite{matsumura}.
A topological analogue of isomorphism (\ref{clin}) for certain linear algebraic groups over $\C$ appears in the unpublished note (\cite{bd}, Theorem 0.11).
\end{rema}

\begin{ex}\rm Consider the special case $G={\rm SL}_2$ and $H=\G$, with $H$ viewed as the diagonal subtorus of $G$. Since $H$ is reductive, the quotient $X=G/H$ is affine of dimension 2 (in fact, it is known to be a quadric in $\A^3$; see e.g. \cite{dukraft}, Example 8.4). By simply connectedness of ${\rm SL}_2$ the complex $\mathcal{C}_{X\lin}$ reduces to $[0\to \Z\oplus \Z\to \Z]$ with the addition of components as second map. Therefore $H_{1}(\mathcal{C}_{X\lin})\cong\Z$ and Corollary \ref{corpigp} gives $\pi_2( X^{\wedge (p')},1)\cong \Z_{(p')}$.
\end{ex}

We thus have an example of a nontrivial second homotopy group for the $p'$-completed \'etale homotopy type of an affine scheme. Note, however, that due to a general theorem of Achinger \cite{achinger} the \'etale homotopy groups of a connected affine scheme of characteristic $p>0$ always vanish in degrees $>1$. Therefore taking the $p'$-completion of the \'etale homotopy type is crucial if one is aiming at results in positive characteristic that are in accordance with those in characteristic 0. The subtlety of having to complete the \'etale homotopy type {\em before} taking homotopy groups is a phenomenon that only occurs for higher homotopy groups because for $i=1$ one has $\pi_1( X^{\wedge (p')})\cong  \pi_1(X)^{(p')}$ as already recalled above. Thus the correct higher analogues of the  prime-to-$p$ fundamental group are the groups $\pi_i( X^{\wedge (p')})$.

We are grateful to Piotr Achinger, Mattia Talpo and Burt Totaro for very helpful exchanges.

\section{Facts from \'etale homotopy theory}\label{secethomfacts}

In this section we collect facts from \'etale homotopy theory to be used in the proofs below. From now on we shall be unforgivably sloppy in notation and shall ignore base points. The notation $X_{\rm et}$ will stand for the \'etale homotopy type of a scheme $X$; when completing it we shall drop the subscript `et'. The $i$-th homotopy group of $X_{\rm et}$ is the $i$-th \'etale homotopy group $\pi_i(X)$ of $X$.

\begin{facts}\label{ethomfacts}\rm ${}$\medskip

\noindent (1) If $X$ is a connected Noetherian normal (more generally, geometrically unibranch) scheme, the \'etale homotopy type is profinite, and therefore the groups $\pi_i(X)$ are profinite groups for $i>0$. \medskip

\noindent (2) A finite surjective radicial morphism $X \to Y$ of Noetherian schemes induces an isomorphism of \'etale homotopy types, hence of \'etale homotopy groups. \medskip

\noindent (3) If $S$ is the spectrum of a discrete valuation ring with separably closed residue field $k$ and fraction field $K$, then for a {\em smooth proper} $S$-scheme $X$ with connected fibres there exists a specialization map
$
X_{K,\rm et}\to X_{k, \rm et}
$
of \'etale homotopy types, inducing isomorphisms
$$
\pi_i^{}(X_K^{\wedge (p')})\stackrel\sim\to \pi_i(X_k^{\wedge (p')})
$$
where $p={\rm char}(k)$.\medskip

\noindent (4) Let $X\to Y$ be a smooth morphism of normal schemes {\em such that Zariski locally $X$ has a smooth relative normal crossing compactification over $Y$} (see \cite{friedlander}, Definition 11.4 for the precise notion). Given a geometric point $\bar y$ of $Y$ there is a long exact homotopy sequence
$$
\cdots\to\pi_i(X_{\bar y}^\wedge)\to \pi_i(X^\wedge)\to \pi_i(Y^\wedge)\to \pi_{i-1}(X_{\bar y}^\wedge)\to\cdots
$$
where $\wedge$ means completion away from the residue characteristics of $Y$.\medskip

\noindent (5) A finite \'etale cover $\phi:\, X\to Y$ of schemes induces isomorphisms $\pi_i(X)\stackrel\sim\to \pi_i(Y)$ for $i\geq 2$. If moreover $X$ and $Y$ are normal schemes of exponential characteristic $p\geq 0$ and the degree of $\phi$ is prime to $p$, then also $\pi_i(X^{\wedge (p')})\stackrel\sim\to \pi_i(Y^{\wedge (p')})$ for $i\geq 2$.
\end{facts}

Here (1) is proven in \cite{artinmazur}, \S 11, (2) is a consequence of \cite{sga1}, IX, Theorem 4.10,(3) is \cite{artinmazur}, Corollary 12.13, and (4) appears in \cite{frmanu} and  \cite{friedlander}, \S 11 (see also the proof of \cite{ss} Proposition 2.8). Finally, (5) in the uncompleted case follows in view of \cite{ss}, Lemma 2.1 from the analogous statement for topological covers, and the $p'$-completed case results from applying \cite{artinmazur}, Theorem 4.11.

We shall also need an invariance property under base field extensions (in fact, only the characteristic 0 case will be used).

\begin{prop}\label{invar} If $K\supset k$ is an extension of separably closed fields of characteristic $p\geq 0$ and $X$ is a connected scheme of finite type over $k$, the morphism $X_{K,\rm et}^{\wedge (p')}\to X_{\rm et}^{\wedge (p')}$ of $p'$-completed \'etale homotopy types is an isomorphism and hence the natural maps
$$
\pi_i(X_K^{\wedge (p')})\to \pi_i(X^{\wedge (p')})
$$
are isomorphisms for all $i \geq 0$.\end{prop}

\begin{proof} For $X$ proper  this is proven in \cite{artinmazur}; in fact in that case it is enough to consider profinite completions instead of $p'$-completions. The same argument works under our assumptions, using as geometric inputs (\cite{sga4}, expos\'e XVI, Corollary 1.6) instead of an application of the proper base change theorem in \'etale cohomology as well as the case $i=1$ treated in (\cite{org}, Corollary 4.5). (To be honest, the statement in \cite{org} assumes that $K$ and $k$ are algebraically closed but the result holds in the separably closed case as well thanks to Fact \ref{ethomfacts} (2).)
\end{proof}

Finally we recall the following basic property of curves.

\begin{prop}\label{curv}
If $X$ is a smooth connected affine curve over an algebraically closed field of characteristic $p \geq 0$, then $\pi_i(X)=0$ and $\pi_i(X^{\wedge (p')})=0$ for $i\geq 2$.\end{prop}

\begin{proof}
The uncompleted case is Proposition 15 in \cite{schmidt}. In the $p'$-completed case we may proceed similarly, by passing to the maximal prime-to-$p$ pro-\'etale cover $\widetilde X\to X$. By Fact \ref{ethomfacts} (5) the natural \hbox{maps} $\pi_i(\widetilde X^{\wedge (p')})\to \pi_i(X^{\wedge (p')})$ are isomorphisms for $i\geq 2$; moreover, the groups $\pi_i(\widetilde X^{\wedge (p')})$ are trivial for $i=0,1$. To prove that they are trivial also for $i>1$, we may apply (\cite{artinmazur}, Theorem 4.3) which reduces the statement to the classical fact that the \'etale cohomological dimension of smooth affine curves is 1.
\end{proof}

\section{Geometric fibrations and homogeneous spaces}

We now establish a number of cases where the geometric assumption of Fact \ref{ethomfacts} (4) can be verified directly in the context of homogeneous spaces and therefore the fibration sequence exists. These cases will be used in subsequent sections for proving Theorems \ref{pi2p} and \ref{homex} of the introduction.

The key lemma is the following.

\begin{lem}\label{lem adjoint}
Let $k$ be algebraically closed, and $f : Y \to X$ a $G$-equivariant morphism of left homogeneous spaces of an algebraic group $G$ such that  moreover $f$ is a right $X$-torsor under a connected $k$-group $H$. Assume that every right $H$-torsor over a field extension $L\supset k$ has a smooth normal crossing compactification over $L$.

Then Zariski locally $Y$ has a smooth relative normal crossing compactification over $X$. Consequently, there is a long exact sequence
$$\cdots \to  \pi_{i+1}( X^{\wedge (p')}) \to \pi_i( H^{\wedge (p')}) \to \pi_i( Y^{\wedge (p')}) \to \pi_i( X^{\wedge (p')}) \to \cdots $$
of \'etale homotopy groups.
\end{lem}

\begin{dem} In view of Fact \ref{ethomfacts} (4) it suffices to prove the first statement.
Let $\eta \in X$ be the generic point. The generic fiber $Y_\eta \to \spec K(X)$ is a right $K(X)$-torsor under $H$, and therefore by assumption there exists a smooth normal crossing compactification $\iota : Y_\eta \to Y_\eta^c$ over $\spec K(X)$. It follows that there is a nonempty Zariski open subset $U \subset X$ such that the restriction $f_U : Y_U \to U$ of $f$ has a smooth relative normal crossing compactification $Y_U^c$ over $U$.
For $g \in G(k)$ consider the translate $gU\subset X$ of $U$ in $X$. Multiplication by $g$ induces an isomorphism $U\stackrel\sim\to gU$ of open sets in $X$, whence also an isomorphism between the $H$-torsors $f_U:\,Y_U\to U$ and $f_{gU}:\,Y_{gU}\to gU$. Since $f_U$ has a smooth relative normal crossing compactification over $U$, we obtain one for $f_{gU}$ by transport of structure. Finally, by transitivity of the $G$-action on $X$ every closed point of $X$ is contained in some $gU$. As $X$ is a finite-dimensional noetherian scheme, this shows that the union of the open sets $gU$ for all $g\in G(k)$ is the whole of $X$.
\end{dem}

Now we collect cases where the geometric condition of the lemma imposed on $H$ is satisfied.

\begin{exs} \label{ex solvable sab} ${}$ \rm
\begin{enumerate}
  \item Assume $H$ is linear, connected and solvable. Then $H$ is special in the sense of Serre (i.e. $H$-torsors are Zariski locally trivial; see \cite{efa}, \S 4.4, Proposition 14). So we only have to find a smooth normal crossing compactification of $H$ over $k$, which exists since as a variety it is isomorphic to the product of a torus and an affine space.
  \item If $H$ is a semi-abelian variety, then writing $H$ as an extension of an abelian variety by a torus and considering the projective bundle associated to the toric bundle we see that the condition of the lemma is satisfied.
\end{enumerate}
\end{exs}

Another case is contained in:

\begin{lem} \label{lem torsor compact}
Let $H$ be an adjoint $k$-group, $L|k$ a field extension and $Y$ a (right) $L$-torsor under $H$.
Then $Y$ has a smooth normal crossing compactification over $L$.
\end{lem}

\begin{dem}
The group $H$ admits a wonderful compactification ${\iota:H\hookrightarrow H^c}$ over $k$ (see \cite{BK}, \S 6.1). In particular, $H^c$ is a smooth projective $k$-variety containing $H$ as the complement of a normal crossing divisor, and the right action (by multiplication) of $H$ on itself extends to $H^c$.
 Denote by $L^s$ a separable closure of $L$, and pick $y_0 \in Y(L^s)$. The point $y_0$ defines a natural isomorphism of $L^s$-varieties $\varphi_0 : H_{L^s} \xrightarrow{\sim} Y_{L^s}$. Consider the open embedding $Y_{L^s} \to {H}_{L^s}^c$ of $L^s$-varieties defined by  $\iota_0 := \iota \circ \varphi^{-1}_0$.

Since $Y$ is an $L$-torsor under $H$, for all $\gamma \in \gal(L^s|L)$ there exists a unique $h_\gamma \in H(L^s)$ such that ${\gamma (y_0)} = y_0 \cdot h_\gamma$. We now twist the $\gal(L^s|L)$-action on ${H}_{L^s}^c$ by the cocycle $\gamma \mapsto h_\gamma$, i.e. we make ${\gamma \in \gal(L^s|L)}$ act on $x \in  {H^c}(L^s)$ by $x \mapsto h_\gamma \cdot {\gamma (x)}$.
Galois descent (see for instance \cite{BLR}, \S 6.2, Example B) implies that ${H}_{L^s}^c$ equipped with its twisted Galois action descends to a smooth projective $L$-variety $Y^c$, and the $\gal(L^s|L)$-equivariant morphism $\iota_0 : Y_{L^s} \to H^c_{L^s}$ comes from an $L$-morphism $\iota_Y : Y \to Y^c$. The morphism $\iota_Y$ is an open immersion and the complement is a normal crossing divisor, since the normal crossing property is local for the \'etale topology.
\end{dem}

\section{Proof of Theorem \ref{pi2p}}

In this section we prove Theorem \ref{pi2p}. As a warm-up, we begin with:

\begin{lem}\label{cs} If $G$ is a connected solvable linear algebraic group over an algebraically closed field of characteristic $p\geq 0$, then $\pi_i( G^{\wedge (p')})=0$ for $i>1$. If moreover $G$ is unipotent, we also have  $\pi_1^{(p')}( G)=0$.
\end{lem}

\begin{dem}
By Proposition \ref{curv} we have $\pi_i(\G^{\wedge (p')})=\pi_i(\Ga^{\wedge (p')})=0$ for $i>1$. Also, $\pi_1^{(p')}(\Ga)=0$ as is well known. Thus by successive application of the fibration sequence of Fact \ref{ethomfacts} (4) for $i>1$ we get $\pi_i( T^{\wedge (p')})=0$ for a torus $T$ and for $i\geq 1$ we get $\pi_i(({\A^n})^{\wedge (p')})=0$ for affine $n$-space $\A^n$. Now the underlying $k$-variety of a unipotent $G$ is just an affine $n$-space, so the second statement follows. For the first, note that a connected $G$ is isomorphic as a variety to the direct product of a torus and an affine space, so we conclude by another application of the fibration sequence.
\end{dem}

Next we consider the case of linear algebraic groups over $\C$. For this we need a comparison result for classical and \'etale homotopy groups. In the statement below, the notation $\pi_i^{\rm top}$ stands for the $i$-th classical homotopy group of a topological space.

\begin{prop}\label{compth} Let $G$ be a connected smooth algebraic group over $\C$. For all $i>0$ there are natural maps $$\pi_i^{\rm top}(G(\C))\to \pi_i(G)$$ inducing isomorphisms $$\pi_i^{\rm top}(G(\C))^\wedge\stackrel\sim\to \pi_i(G) = \pi_i(G^\wedge),$$ where $\wedge$ denotes profinite completion.
\end{prop}

We thank Burt Totaro for his help with the proof below.\medskip

\begin{dem}
By (\cite{artinmazur}, Theorem 12.9 and Corollary 12.10) for any geometrically unibranch connected normal scheme $X$ of finite type over $\C$  there is a comparison map $X_{\rm cl}\to X_{\rm et}$ of homotopy types inducing an isomorphism of the profinite completion $X_{\rm cl}^\wedge$ of $X_{\rm cl}$ with $X_{\rm et} = X^\wedge$. Here $X_{\rm cl}$ computes the classical homotopy groups of $X(\C)$, i.e. $\pi_i(X_{\rm cl})=\pi_i(X(\C))$.

It remains to see that for $X=G$ the natural maps $\pi_i(G_{\rm cl})^\wedge\to \pi_i(G_{\rm cl}^\wedge)$ are isomorphisms for $i>0$. To see this, recall first that the fundamental group of $G(\C)$ is abelian (this is true for every topological group), and moreover all of its homotopy groups are finitely generated abelian groups. Indeed, the integral homology groups of $G(\C)$ are finitely generated  (since so are those of a maximal compact subgroup which is a deformation retract of $G(\C)$). On the other hand, again since $G(\C)$ is a topological group, it is a nilpotent space in the sense of homotopy theory (see e.g. \cite{mayponto}, Corollary 1.4.5 and Definition 3.1.4). Therefore its homotopy groups are also finitely generated (see e.g. \cite{mayponto}, Theorem 4.5.2). Now our claim about completions of homotopy groups follows from (\cite{sullivan}, Theorem 3.1).
\end{dem}

\begin{cor}\label{corc}
If $G$ is a connected linear algebraic group over $\C$, then $\pi_2(G)=0$.
\end{cor}

\begin{proof}
By the proposition we are reduced to proving $\pi_2^{\rm top}(G(\C))=0$. At this point we invoke Cartan's theorem: we have $\pi_2^{\rm top}(K)=0$ for a maximal compact subgroup $K$ in the underlying real Lie group of $G(\C)$. But $K$ is a deformation retract of $G(\C)$, whence the corollary.
\end{proof}

Now we can treat the case of linear groups in general.

\begin{prop}\label{lin} If $G$ is a connected linear algebraic group over an algebraically closed field of characteristic $p\geq 0$, then $\pi_2( G^{\wedge (p')})=0$.
\end{prop}

\begin{dem}
Let $G_u$ be the unipotent radical of $G$. Using Example \ref{ex solvable sab}, we may apply Lemma \ref{lem adjoint} to the $G_u$-torsor over $G/G_u$ defined by the extension
$$
1\to G_u\to G\to G/G_u\to 1
$$
and consider the associated homotopy sequence. Since  $\pi_2( G_u^{\wedge (p')})=0$ by Lemma \ref{cs}, we see that in order to prove $\pi_2( G^{\wedge (p')})=0$ we may replace $G$ by $G/G_u$ and hence assume from now on that $G$ is reductive.

In characteristic 0 Proposition \ref{invar} allows us to reduce to the case $k=\C$ which is contained in  Corollary \ref{corc}.
To treat the case $p>0$, recall that $G$ extends to a reductive group scheme $\widetilde G$ over the Witt ring $W(k)$ by (\cite{sga3}, Expos\'e XXV, Corollaire 1.3). As $W(k)$ is strictly henselian, there exists a Borel subgroup $\widetilde B\subset \widetilde G$ by (\cite{sga3}, Expos\'e XXII, Corollaire 5.8.3 (i)) and we may consider the quotient $\widetilde G/\widetilde B$. Denote the geometric generic fibres of $\widetilde G$ and $\widetilde B$ by $G_0$ and $B_0$, respectively, and let $B$ be the special fibre of $\widetilde B$. Writing $G$ as a $B$-torsor over $G/B$ we can again conclude from Lemma \ref{lem adjoint} that the quotient map $G\to G/B$ sits in a long exact fibration sequence (see Example \ref{ex solvable sab}). The same is true for the map $G_0\to G_0/B_0$, whence the horizontal maps in the exact commutative diagram
$$
\begin{CD}
0 @>>> \pi_2( G^{\wedge (p')}) @>>> \pi_2({(G/B)}^{\wedge (p')}) @>>> \pi_1(B)^{(p')}\\
&& && @AA{\cong}A @AA{\cong}A \\
0 @>>> \pi_2( G_0^{\wedge (p')}) @>>> \pi_2({(G_0/B_0)}^{\wedge (p')}) @>>> \pi_1(B_0)^{(p')}.
\end{CD}
$$
The zeros on the left come from the vanishing of $\pi_2( B^{\wedge (p')})$ and $\pi_2( B_0^{\wedge (p')})$ implied by  Lemma \ref{cs}. The middle vertical isomorphism is that of Fact \ref{ethomfacts} (3) applied to the proper smooth $W(k)$-scheme $\widetilde G/\widetilde B$. Given compatible maximal tori $T\subset B$ and $T_0\subset B_0$ we have isomorphisms $\pi_1(B)^{(p')}\cong \pi_1(T)^{(p')}$ and $\pi_1(B_0)^{(p')}\cong \pi_1(T_0)^{(p')}$ since $B$ is a product of $T$ with some affine space and similarly for $B_0$. (Recall that $\pi_1^{(p')}$ is compatible with direct products \cite{org} and $\pi_1^{(p')}(\A^1)=0$.) The right vertical map is therefore identified with the specialization map $\pi_1(T_0)^{(p')}\to \pi_1(T)^{(p')}$ coming from the specialization theory of the tame fundamental group of split tori and is an isomorphism by compatibility with products and the case of $\G$ (see \cite{ovi}, Th\'eor\`eme 4.4). Now the diagram implies the existence of an isomorphism $\pi_2( G_0^{\wedge (p')})\stackrel\sim\to \pi_2( G^{\wedge (p')})$, but $\pi_2( G_0^{\wedge (p')})=0$ by the characteristic 0 case.
\end{dem}

\begin{rema}\rm With the above notation, the proof shows that for $G$ reductive we in fact have isomorphisms $\pi_i( G_0^{\wedge (p')})\stackrel\sim\to \pi_i( G^{\wedge (p')})$ for all $i$, i.e. the $p'$-completed \'etale homotopy types of $G$ and $G_0$ are weakly equivalent. This follows from continuing the fibration sequence: for $i>2$ the argument is straightforward and for $i=1$ one has to use the fact that flag varieties are simply connected over $\C$, hence over a field of characteristic 0. For a related statement, see (\cite{frK}, Proposition 2.8).
\end{rema}

\begin{lem}\label{av}
If $A$ is an abelian variety over an algebraically closed field $k$ of characteristic $p$, then $\pi_i( A^{\wedge (p')})=0$ for $i>1$.
\end{lem}

\begin{dem}
As in the proof of the previous proposition, in the case $p=0$ we reduce to verifying the claim for $k=\C$ and the usual homotopy groups of $A(\C)$. These are trivial by the usual long exact homotopy sequence as then $A(\C)$ is the quotient of the contractible space $\C^g$ by a discrete subgroup. In the case $p>0$ we again proceed by specialization: abelian varieties lift to characteristic 0 (see e.g. \cite{no}, Corollary 3.2) and Fact \ref{ethomfacts} (3) applies.
\end{dem}

\begin{cor}\label{sab}
If $G$ is a semi-abelian variety over an algebraically closed field $k$ of characteristic $p$, then $\pi_i( G^{\wedge (p')})=0$ for $i>1$.
\end{cor}

\begin{dem}
Write $G$ as an extension of an abelian variety $A$ by a torus $T$. The quotient map $G\to A$ satisfies the assumption in Fact \ref{ethomfacts} (4) (consider $G$ as a torus bundle over $A$ and take the associated projective bundle: see Example \ref{ex solvable sab}) and therefore the homotopy exact sequence may be applied to reduce the corollary to Lemmas \ref{cs} and \ref{av}.
\end{dem}

\noindent {\em Proof of Theorem \ref{pi2p}.} 
Observing that the \'etale homotopy type is unaffected by purely inseparable base change (see Fact \ref{ethomfacts} (2)), we may assume $k$ algebraically closed. Then the proof for $p=0$ is the same as in the linear case treated in Proposition \ref{lin}, so assume $p>0$. In this case the largest anti-affine subgroup $G_{\rm ant}$ of $G$ (i.e. the largest closed subgroup $H\subset G$ with $\calo(H)=k$) is a semi-abelian variety central in $G$ by (\cite{aa}, Proposition 2.2). By Lemma \ref{lem adjoint} and Example \ref{ex solvable sab}, we therefore have a fibration sequence for the quotient map $G\to G/G_{\rm ant}$. Moreover, the quotient $G/G_{\rm ant}$ is linear (see e.g. \cite{brionbook}, Theorem 3.2.1), and therefore the theorem follows from Proposition \ref{lin} and Corollary \ref{sab}.\enddem

\begin{rema}\label{gantrema}\rm By a theorem of Rosenlicht (see e.g. \cite{brionbook}, Theorem 1.2.1) the subgroup $G_{\rm ant}\subset G$ used in the above proof is the smallest normal subgroup $H\subset G$ such that $G/H$ is affine. If $G^{\rm lin}$ has trivial unipotent radical, then so does $G/G_{\rm ant}$, and the inverse image in $G$ of a maximal torus of $G/G_{\rm ant}$ defines a maximal semi-abelian subvariety ${\rm SA}_G\subset G$ as considered before Proposition \ref{proppi1G} (see \cite{dem:brauer}, \S 4.1). It is not a normal subgroup of $G$ in general.
\end{rema}

\section{The homotopy exact sequence down to degree 3}

In this section we establish the homotopy exact sequence of Theorem \ref{homex} in degrees $\geq 3$ and prove some auxiliary statements that will also serve in the low-degree part. We shall assume throughout that the base field $k$ is algebraically closed and that the connected subgroup $H\subset G$ is smooth, which is allowed by Fact \ref{ethomfacts}(2) and Remark \ref{Hred}.

Recall that we have to construct the fibration sequence for the quotient map $G\to X$ with stabilizer $H$. The proof will proceed by breaking up the $H$-torsor $G\to X$ in pieces. To this end, let us introduce some notation. Denote by $H\u$ the unipotent radical of $H\lin$, by ${H\red := H\lin / H\u}$ the reductive quotient of $H\lin$ and by $H\ab$ the quotient abelian variety $H/H\lin$.

\begin{lem}\label{devis} Let $G$ be a connected algebraic $k$-group, with a smooth connected closed subgroup $H\subset G$, and let $X:=G/H$. Given a Borel subgroup $B \subset H\red$, the quotient map $G\to X$ factors as a sequence
\[
\xymatrix{
G \ar[r]^{H\u} & W \ar[rd]^B \ar[rr]^{H\red} & & Y \ar[r]^{H\ab} & X \, , \\
&  & Z \ar[ru]^\pi & &
}
\]
where each map labelled by a $k$-group is a (right) torsor under this group, while the morphism $\pi : Z \to Y$ is smooth and proper.
\end{lem}

\begin{dem}
Set $W := G/H\u$, $Y := G/H\lin$ and $Z := W/B$. The only property that requires a proof is that of the morphism $\pi$, which is a consequence of \cite{sga3}, expos\'e XXII, Corollary 5.8.3, and of the fact that properness and smoothness for a morphism can be checked \'etale-locally.
\end{dem}

\begin{lem}\label{pipi}
With notation as in the lemma above, there are isomorphisms
\begin{enumerate}
\item $\pi_i(G^{\wedge (p')}) \xrightarrow{\sim} \pi_i(Z^{\wedge (p')})$ for all $i\geq 3$;
\item $\pi_i(Y^{\wedge (p')}) \xrightarrow{\sim} \pi_i(X^{\wedge (p')})$ for all $i\geq 3$.
\end{enumerate}
\end{lem}

\begin{dem}
\begin{enumerate}
\item The morphism $G \to Z$ is the composition of two morphisms that are torsors under solvable groups ($H\u$ and $B$), hence by lemmas \ref{cs}, \ref{lem adjoint} and Example \ref{ex solvable sab}, we get the isomorphisms.
\item The map $Y \to X$ is a torsor under an abelian variety, so as in the proof of Theorem \ref{pi2p}, we may apply Lemma \ref{lem adjoint} and Example \ref{ex solvable sab} to derive isomorphism (2) from Lemma \ref{av}.
\end{enumerate}
\end{dem}

Now we can prove as promised:

\begin{prop}
With notation as in Lemma \ref{devis} there is a long exact sequence of the shape
$$
\cdots \to  \pi_{4}( X^{\wedge (p')}) \to \pi_3( H^{\wedge (p')}) \to \pi_3( G^{\wedge (p')}) \to \pi_3( X^{\wedge (p')}) \to 0.
$$
\end{prop}

\begin{dem}
By Fact 2.1.(5), the smooth proper morphism $\pi : Z \to Y$ gives rise to a long exact fibration sequence
\begin{equation}\label{exact sequence Had}\cdots \to \pi_{i+1}( Y^{\wedge (p')}) \to \pi_i\left(\left(H/B\right)^{\wedge (p')}\right) \to \pi_i( Z^{\wedge (p')}) \to \pi_i( Y^{\wedge (p')}) \to \cdots \end{equation}
which may be rewritten for $i\geq 3$, using Lemma \ref{pipi}, as
\begin{equation}\label{pipieq}\cdots\to \pi_{i+1}( X^{\wedge (p')}) \to \pi_i\left(\left(H/B\right)^{\wedge (p')}\right) \to \pi_i( G^{\wedge (p')}) \to \pi_i( X^{\wedge (p')}) \to \cdots \end{equation}
Since $H \to H/B$ is a $B$-torsor, Lemmas \ref{cs}, \ref{lem adjoint} and Example \ref{ex solvable sab} provide isomorphisms $\pi_i( H^{\wedge (p')})\stackrel\sim\to\pi_i\left(\left(H/B\right)^{\wedge (p')}\right)$ for $i \geq 3$, so
we may replace $\pi_i\left(\left(H/B\right)^{\wedge (p')}\right)$ by $\pi_i( H^{\wedge (p')})$ in (\ref{pipieq}).
\end{dem}

\section{The fibration sequence in low degree}\label{secsec}

In this section we establish the remaining statements announced in the introduction. We assume $k$ is algebraically closed and introduce the following classical notation. Let $H$ be a smooth connected algebraic $k$-group. Let $H\ss \subset H\red$ be the derived subgroup of $H\red$, which is semisimple. Consider $H\sab := (H^0 / H\u) / H\ss$, the maximal semi-abelian quotient of $H / H\u$. It is isomorphic to the quotient of the semi-abelian variety ${\rm SA}_{H}\subset H / H\u$ introduced before Proposition \ref{proppi1G} by a maximal torus in $H\ss$.
Finally, let $Z_{H\ss} \subset H\ss$ be the center of $H\ss$ and let $H\ad := H\ss / Z_{H\ss}$ be the adjoint quotient of $H\ss$.

\begin{lem}\label{devis} For a smooth connected algebraic $k$-group $G$, a closed smooth connected $k$-subgroup $H\subset G$ and $X:=G/H$, the quotient map $G\to X$ factors as a sequence
$$G \xrightarrow{H\u} W \xrightarrow{Z_{H\ss}} V \xrightarrow{H\ad} Y \xrightarrow{H\sab} X \, ,$$
where each map is labelled by a $k$-group under which it a (right) torsor.
\end{lem}

\begin{dem}
Set $W := G/H\u$. To define $V$, denote by $Z\u_{H\ss}$ the inverse image of $Z_{H\ss}$ in $H\lin$, and set $V := G/Z\u_{H\ss}$. Finally, set $Y := G / H\ssu$, where $H\ssu$ is the inverse image of $H\ss$ in $H\lin$. Since by construction $W\to Y$ is an $H\ss$-torsor, we conclude that $V\to Y$ is indeed an $H\ad$-torsor.
\end{dem}

With notation as in Lemma \ref{devis}, set $X':=W/{\rm SA}_{H}$, where ${\rm SA}_{H}\subset H/H^{\rm u}$ is the maximal semi-abelian subvariety introduced before Proposition \ref{proppi1G}. On the other hand, recall that by construction we have $$X \cong (G/H^{\rm u})/(H/H^{\rm u})=W/(H/H^{\rm u}),$$ so there is a natural map $X'\to X$.

\begin{lem}\label{lemY'}
The map $X'\to X$ induces a canonical isomorphism
$$
\pi_1(X')^{(p')}\stackrel\sim\to  \pi_1(X)^{(p')}
$$
and a canonical exact sequence
\begin{equation}
0\to {T_{H\sc}}_*\otimes\Z_{(p')}(1)\to  \pi_2( X'^{\wedge (p')})\to  \pi_2( X^{\wedge (p')})\to 0.
\end{equation}
\end{lem}

\begin{dem}
Consider first the quotient $Y' := W / T_{H\ss}\cong V/T_{H\ad}$, where $T_{H\ss}\subset H\ss$ and $T_{H\ad}\subset H\ad$ are compatible maximal tori. The map $V\to Y$ factors through $Y'$ and $V\to Y'$ is a torsor under $T_{H\ad}$. Using Lemma \ref{lem adjoint}, Example \ref{ex solvable sab} (for the first line) and Lemma \ref{lem torsor compact} (for the second line), we get the following commutative exact diagram
\begin{displaymath}
 \xymatrix{
0 \ar[r] & \pi_2(V^{\wedge (p')}) \ar[d]^= \ar[r]& \pi_2( Y'^{\wedge (p')}) \ar[r] \ar[d] & \pi_1(T_{H\ad})^{(p')} \ar[r] \ar[d] &  \pi_1(V)^{(p')} \ar[r] \ar[d]^= &  \pi_1(Y')^{(p')} \ar[r] \ar[d] & 0 \\
0 \ar[r] & \pi_2(V^{\wedge (p')}) \ar[r] &\pi_2( Y^{\wedge (p')}) \ar[r] &  \pi_1(H\ad)^{ (p')} \ar[r] &  \pi_1(V)^{(p')} \ar[r] &  \pi_1(Y)^{(p')} \ar[r] & 0
 }
\end{displaymath} where the zeros on the left come from Proposition \ref{lin}. By (\cite{dem}, Proposition 3.10), the third vertical map is surjective with kernel $ \pi_1(T_{H\sc})^{(p')}\cong { {T_{H\sc}}_*\otimes\Z_{(p')}(1)}$. So a diagram chase gives an exact sequence
\begin{equation}\label{ex1}
0\to {T_{H\sc}}_*\otimes\Z_{(p')}(1)\to  \pi_2( Y'^{\wedge (p')})\to  \pi_2( Y^{\wedge (p')})\to 0
\end{equation}
as well as an isomorphism
\begin{equation}\label{iso}
\pi_1(Y')^{(p')}\to \pi_1(Y)^{(p')}.
\end{equation}
Now the quotient $X' = W / {\rm SA}_{H}$ gives rise to a right torsor $Y'\to X'$ under the semi-abelian variety $H\sab$. Since on the other hand $Y \cong W/H\ss$ and $Y\cong W/(H/H\u)$, we have a commutative diagram
$$
\begin{CD}
Y' @>{H\sab}>> X'\\
@VVV @VVV \\
Y @>{H\sab}>> X
\end{CD}
$$
of right torsors under $H\sab$. The associated homotopy exact sequences constructed using Lemma \ref{lem adjoint} and Example \ref{ex solvable sab} give rise to a commutative exact diagram (see Corollary \ref{sab})
{\tiny $$
\xymatrix{
0 \ar[r] & \pi_2(Y'^{\wedge (p')}) \ar[r] \ar[d] & \pi_2(X'^{\wedge (p')}) \ar[r] \ar[d] & \pi_1(H\sab)^{(p')} \ar[r] \ar[d]^= & \pi_1(Y')^{(p')}  \ar[d]^\cong \ar[r] & \pi_1(X')^{(p')} \ar[d] \ar[r] & 0 \\
0 \ar[r] & \pi_2(Y^{\wedge (p')})  \ar[r] & \pi_2(X^{\wedge (p')}) \ar[r] & \pi_1(H\sab)^{(p')}  \ar[r] & \pi_1(Y)^{(p')} \ar[r] & \pi_1(X)^{(p')} \ar[r] & 0  \,
}
$$}

\noindent where the fourth vertical map is an isomorphism by (\ref{iso}). The lemma follows from the diagram and
exact sequence (\ref{ex1}).\end{dem}

Using the lemma we can already determine $\pi_1(G)^{(p')}$ as announced in the introduction.\medskip

\noindent{\em Proof of Proposition \ref{proppi1G}.}
With notation as in the previous proof, the torsor $W \to X'$  under the semi-abelian variety $\SA_{H}$ gives rise to an exact fibration  sequence
$$
\pi_2( W^{\wedge (p')}) \to \pi_2( X'^{\wedge (p')}) \to \pi_1({\rm SA}_{H})^{(p')} \to \pi_1(W)^{(p')} \to \pi_1(X')^{(p')}
$$
in a by now familiar fashion.
The morphism $G \to W$ is a torsor under the  unipotent group $H\u$, and therefore Lemmas \ref{cs}, \ref{lem adjoint} and Example \ref{ex solvable sab} imples that $\pi_2(W^{\wedge (p')})=0$, hence the previous sequence can be written as:
\begin{equation} \label{ses Y'}
0\to \pi_2( X'^{\wedge (p')}) \to \pi_1({\rm SA}_{H})^{(p')} \to \pi_1(G)^{(p')} \to \pi_1(X')^{(p')}.
\end{equation}
Now set $G=H$. In this case $X$ is a point and the lemma above gives isomorphisms $\pi_1(X')^{(p')}\stackrel\sim\to  \pi_1(X)^{(p')}=0$
and ${T_{G\sc}}_*\otimes\Z_{(p')}(1)\stackrel\sim\to  \pi_2(X'^{\wedge (p')})$. It remains to recall that $\pi_1(\SA_G)^{(p')}\cong T_{(p')}(\SA_G)$. For an abelian variety this well-known fact can be found e.g. in (\cite{szamuely}, Theorem 5.6.10). The semi-abelian case is proven in the same way using (\cite{bsz}, Proposition 1.1).
\enddem

Next we  prove Corollary \ref{corpigp}.

\begin{prop}  There is a canonical isomorphism of abelian profinite  groups
$$\pi_2( X^{\wedge (p')})\cong H_{1}(\mathcal{C}_{X,p'}),  $$
where
$$\mathcal{C}_{X,p'} := [{T_{H\sc}}_* \otimes \Z_{(p')}(1) \to T_{(p')}(\SA_{H}) \oplus ({T_{G\sc}}_* \otimes \Z_{(p')}(1) )\to T_{(p')}(\SA_G)].$$
\end{prop}

\begin{dem}
Substituting the formula of Proposition \ref{proppi1G} in exact sequence (\ref{ses Y'}) gives
$$
0\to \pi_2(X'^{\wedge (p')}) \to T_{(p')}({\rm SA}_{H}) \to \cok({T_{G\sc}}_*\otimes\Z_{(p')}(1)\to T_{(p')}(\SA_G))).
$$
Now apply the exact sequence of Lemma \ref{lemY'} to get the desired formula.

\end{dem}

We finally prove the remaining part of Theorem \ref{homex}, that is:

\begin{prop}\label{homexdown} In the situation of Theorem \ref{homex} there is an exact sequence $$
0 \to \pi_2( X^{\wedge (p')}) \to \pi_1(H)^{(p')} \to \pi_1(G)^{(p')} \to \pi_1(X)^{(p')} \to 0. $$
\end{prop}

\begin{dem}
The above proof of Corollary \ref{corpigp} also yields the exact sequence
$$
0 \to \pi_2( X^{\wedge (p')}) \to \pi_1(H)^{(p')} \to \pi_1(G)^{(p')}
$$
in view of Proposition \ref{proppi1G} applied to $G$ and $H$. So, using Lemma \ref{lemY'}, we may rewrite part of exact sequence (\ref{ses Y'}) as
$$
\pi_1({\rm SA}_{H})^{(p')} \to \pi_1(G)^{(p')} \to \pi_1(X)^{(p')}
$$
which yields the exactness of the sequence of the proposition at $\pi_1(G)^{(p')}$ since the map $\pi_1({\rm SA}_{H})^{(p')} \to \pi_1(G)^{(p')}$ factors through $\pi_1(H)^{(p')}$. Finally, the map $\pi_1(G)^{(p')} \to \pi_1(X)^{(p')}$ is surjective by \cite{bsz}, Theorem 1.2 (a).
\end{dem}

\end{document}